\documentclass{amsart}

\newtheorem{theorem}{Theorem}[section]
\newtheorem{lemma}[theorem]{Lemma}
\newtheorem{corollary}[theorem]{Corollary}

\theoremstyle{definition}

\theoremstyle{remark}
\newtheorem{remark}[theorem]{Remark}

\numberwithin{equation}{section}

\newcommand{\g}{\mathfrak }

\newcommand{\vp}{\varphi}
\newcommand{\m}{\mathcal }
\newcommand{\op}{\operatorname } 

\newcommand{\E}{\operatorname{End}({M})}

\begin{document}

\title{A note on cocycle-conjugate endomorphisms of von Neumann algebras}

\author{Remus Floricel}
\address{Department of Mathematics and Statistics, University of Ottawa, Ottawa, ON, Canada, K1N 6N5}

\email{floricel@matrix.cc.uottawa.ca}

\subjclass{Primary 46L10, 46L40}

\keywords{von Neumann algebras, endomorphisms, cocycle-conjugacy, shifts}

\begin{abstract}We show that two cocycle-conjugate endomorphisms of an arbitrary von Neumann algebra that satisfy certain stability conditions are conjugate endomorphisms, when restricted to some specific von Neumann subalgebras. As a consequence of this result, we obtain a new criterion for conjugacy of Powers shift endomorphisms acting on factors of type $\rm{I}\sb{\infty}.$
\end{abstract}

\maketitle

\section{Introduction}Let $M$, $N$ be von Neumann algebras. A pair of unital normal *-endomorphisms $\rho$ of $M$ and  $\sigma$ of $N$ are said to be conjugate, if there exists a *-isomorphism $\theta$ from $N$ onto $M$ such that $\rho\circ\theta=\theta\circ\sigma$. They are called cocycle-conjugate endomorphisms, if there exists a unitary element $u$ of $M$ such that $\sigma$ and $\op{Ad}(u)\circ\rho$ are conjugate endomorphisms.\\ 
\indent The study of conjugacy and cocycle-conjugacy classes of (shift-type) endomorphisms acting on factors of type $\rm{I}\sb{\infty}$ and type $\rm{II}\sb{1}$ was initiated by R. Powers in \cite{P88}, and since then, many important classification results have been obtained (see \cite{B96}, \cite{By90}, \cite{L93}, \cite{Pr98} and the references therein). In general, it is difficult to determine whether a pair of endomorphisms are conjugate or cocycle-conjugate, even in the case they are automorphisms. Furthermore, it is often more difficult to find conjugacy invariants than to find cocycle-conjugacy invariants. For example, although a cocycle-conjugacy invariant for endomorphisms acting on a type $\rm{I}\sb{\infty}$ factor can be easily described, the problem of finding  conjugacy invariants is difficult, even for particular classes of endomorhisms (\cite{B96}, \cite{L93}).\\
\indent Our purpose, in this paper, is to find conditions under which two cocycle-conjugate endomorphisms are also conjugate. In fact, we show that two cocycle-conjugate endomorphisms acting on a von Neumann algebra that satisfy certain stability conditions have the property that they are conjugate endomorphisms, when restricted to some specific von Neumann subalgebras. The consequences of this result are two-fold. First of all, we obtain a new criterion for conjugacy of Powers shifts acting on a type $\rm{I}\sb{\infty}$ factor which improves the original criterion of Powers in \cite{P88}. Secondly, we obtain a generalization of an important results of Bures and Yin in \cite{By90}, from the case of a type $\rm{II}\sb{1}$ factor to the case of a von Neumann algebra with a faithful normal tracial state. \\
\indent We close our introduction with a few remarks on notation, most of which is standard. For a von Neumann algebra $M$ acting on a separable Hilbert space $\m{H}$, we denote by $\op{End}(M)$ the semigroup of all injective unital *-endomorphisms of $M$. We note that such endomorphisms are automatically normal. For $\rho\in\op{End}(M)$, we denote by $M_\rho$  the von Neumann algebra generated by $\{\rho^n(M)'\cap M\,:\,n\geq 0\}$. Since \begin{eqnarray}\label{rel1}\rho\left( \rho^n(M)'\cap M\right)\subset  \rho^{n+1}(M)'\cap M,\end{eqnarray} for every positive integer $n$, it follows that the restriction  $\rho|_{M_\rho}$ is an endomorphism of the von Neumann algebra $M_\rho$.\\
Finally, if $\vp$ is a state of $M$, then $M^\vp$ will denote the centralizer of $\vp$, that is the set of all elements $x$ of $M$ such that $\vp(xy)=\vp(yx)$, for all $y\in M$.

\section{The Main Result}
We begin this section with the following lemma which will be used in the proof of the next theorem:
\begin{lemma}\label{ccf1}Let $M$ be a von Neumann algebra and $\rho\in\E$. Assume that there exists a faithful normal state $\vp$ of $M$ which is is $\rho$-invariant {\em (}i.e. $\varphi(\rho(x))=\varphi(x)$ for all $x\in M${\em )}. Then the restriction $\rho|_{M^\vp}$ is an endomorphism of the von Neumann algebra $M^\vp.$
\end{lemma}
\begin{proof}
 Let $\{\sigma _t^\vp\}_{t\in\mathbb{R}}$ be the modular automorphism group of $M$ associated to the state $\vp$. By \cite{S82}, $M^\vp$ is the fixed-point algebra of $\{\sigma _t^\vp\}_{t\in\mathbb{R}}$, and $\sigma _t^{\vp\circ\rho}=\rho^{-1}\circ\sigma _t^\vp\circ\rho,$ for all $t\in\mathbb{R}.$ Thus, if $x\in M^\vp$ and $t\in\mathbb{R}$, we have $$
\sigma_t^\vp(\rho(x))=\rho\circ\sigma_t^{\vp\circ\rho}(x)=\rho\circ\sigma_t^{\vp}(x)=\rho(x),$$ so $\rho(x)\in M^\vp$.
The lemma is proved.\end{proof}

We now state and prove the main result of this paper:
\begin{theorem}\label{cut1}
Let $M$ be a von Neumann algebra, and $\rho,\,\sigma\in\op{End}(M)$ be cocycle-conjugate endomorphisms. Let $\theta$ be a *-automorphism of $M$, and $u\in M$ be a unitary element such that $\theta\circ\sigma=\op{Ad}(u)\circ\rho\circ\theta.$ Assume that there exists a $\rho$-invariant faithful normal state $\vp$ of $M$ such that $u\in M^\vp$. Then the restrictions $\rho\mid_{{M}\sb{\rho}}$ and $\sigma\mid_{{M}\sb{\sigma}}$ are conjugate endomorphisms. 

\end{theorem}
\begin{proof}{(I)} First of all, we assume that $\theta=Id$, so $\sigma=\op{Ad}(u)\circ\rho$. In the first part of the proof, we shall employ an argument from \cite[Theorem 1.2]{By90}. For any positive integer $n$, we define $$
u_n=\prod _{i=1}^n\rho^{i-1}(u).$$
Then $u_n$ is a unitary in $M$, and it is easily seen that it satisfies the following relation:$$
\sigma^n=\op{Ad}(u_n)\circ \rho^n,\; \; \mbox{for all}\;n\in\mathbb{N}.$$
In particular, this relation implies that the mappings $$\op{Ad}(u_n):\rho^n(M)'\cap M\longrightarrow  \sigma^n(M)'\cap M$$ are *-isomorphisms of von Neumann algebras, for each $n\in\mathbb{N}$. On the other hand, for any $x\in\rho^n(M)'\cap M $, we have the following:$$
\op{Ad}(u_{n+1})(x)=\op{Ad}(u_n)\left(\rho^n(u)x\rho^n(u^*)\right)=\op{Ad}(u_n)(x).$$Thus the diagram
$$
\begin{array}{ccc}    \rho^n(M)'\cap M &\stackrel{\op{Ad}(u_n)}{\longrightarrow} &   \sigma^n(M)'\cap M \\
    \cap &&      \cap\\
\rho^{n+1}(M)'\cap M &\stackrel{\op{Ad}(u_{n+1})}{\longrightarrow}  &  \sigma^{n+1}(M)'\cap M \end{array}  $$
 commutes for every $n\in\mathbb{N}$.\\Accordingly, the *-isomorphisms $\op{Ad}(u_n):\rho^n(M)'\cap M\longrightarrow  \sigma^n(M)'\cap M$ can be lifted to a  *-isomorphism $\Theta$ from the $C^*$-algebra $C_\rho$ generated by $\{\rho^n(M)'\cap M\,|\,n\geq 0\}$, onto the $C^*$-algebra $C_\sigma$ generated by $\{\sigma^n(M)'\cap M\,|\,n\geq 0\}$, uniquely determined by the condition\begin{eqnarray}\label{rel2}\Theta\mid_{\rho^n(M)'\cap M}=\op{Ad}(u_n).\end{eqnarray}
\indent Next, we claim that $\Theta$ can be extended to a *-isomorphism from the von Neumann algebra $\overline{C_\rho}^{\,w}=M_\rho$ onto $\overline{C_\sigma}^{\,w}=M_\sigma$, where $\overline{C_\rho}^{\,w}$ denotes the closure of $C_\rho$ in the weak operator topology.
For proving this claim, we consider the GNS-construction $(\pi_\vp,\,\m{H}_\vp,\,\Omega_\vp)$ of $M$ with respect to the faithful normal state $\vp$. As usual, we identify $M$ with $\pi_\vp(M)$, so $\Omega_\vp\in\m{H}_\vp$ is a cyclic and separating vector for $M$, and $\vp(x)=\left<x\Omega_\vp\,,\,\Omega_\vp\right>$, for all $x\in M$.\\
\indent On the other hand, since $u\in M^\vp$, we deduce from Lemma \ref{ccf1} that $$
\vp\circ\op{Ad}(u_n)=\vp,\; \; n\in\mathbb{N}.$$Thus, by norm continuity, \begin{eqnarray}\label{ali1}\vp|_{C_\sigma}\circ\Theta=\vp|_{C_\rho}.\end{eqnarray} Next, we consider the Hilbert subspaces $\m{H}_\rho$ and $\m{H}_\sigma$ of $\m{H}_\vp$, defined by $$
\m{H}_\rho=\overline {\pi_\vp(C_\rho)\Omega_\vp}\;\;\mbox{and}\;\;\m{H}_\sigma=\overline {\pi_\vp(C_\sigma)\Omega_\vp}.$$ Then the $C^*$-algebras $C_\rho$ and $C_\sigma$ act on $\m{H}_\rho$, respectively on $\m{H}_\sigma$, via the representation $\pi_\vp|_{C_\rho}$, respectively $\pi_\vp|_{C_\sigma}$. In fact, up to unitary equivalence, $(\pi_\vp|_{C_\rho}\,,\,\m{H}_\rho)$ and $(\pi_\vp|_{C_\sigma}\,,\,\m{H}_\sigma)$ are exactly the GNS-representations $(\pi_{\vp|_{C_\rho}}\,,\,\m{H}_{\vp|_{C_\rho}})$, respectively $(\pi_{\vp|_{C_\sigma}}\,,\,\m{H}_{\vp|_{C_\sigma}})$ of the $C^*$-algebra $C_\rho$, respectively $C_\sigma$, with respect to the state $\vp|_{C_\rho}$, respectively $\vp|_{C_\sigma}$. \\
\indent By identifying $C_\rho$ with $\pi_{\vp|_{C_\rho}}(C_\rho)$ and $C_\sigma$ with $\pi_{\vp|_{C_\sigma}}(C_\sigma)$, we define an operator $U:\m{H}_\rho\rightarrow \m{H}_\sigma$ by\begin{eqnarray*}
Ux\Omega_\vp=\Theta (x)\Omega_\vp,\; \; \; x\in C_\rho,\end{eqnarray*} and we claim that $U$ is a unitary operator. Indeed, for any $x\in C_\rho$, we have \begin{eqnarray*}
\|Ux\Omega_\vp\|^2&=&\left<Ux\Omega_\vp\,,\, Ux\Omega_\vp\right>\\&=&\left<\Theta(x^*x)\Omega_\vp\,,\,\Omega_\vp\right>\\&\stackrel{(\ref{ali1})}{=}&\left<x^*x\Omega_\vp\,,\,\Omega_\vp\right>\\&=&\|x\Omega_\vp\|^2,\end{eqnarray*}so $U$ is an isometry. But $U$ is surjective, so  $U$ must be a unitary operator.\\
\indent Moreover, for all $x,\,y\in C_\rho$ we have $$
Uxy\Omega_\vp=\Theta(xy)\Omega_\vp=\Theta(x)\Theta(y)\Omega_\vp=\Theta(x)Uy\Omega_\vp,$$
from which we infer that $$
\Theta(x)=UxU^*,\;\;x\in C_\rho.$$ It then follows that $\Theta$ is a weakly continuous mapping, and therefore, it can be extended  by continuity to a *-isomorphism from $M_\rho$ onto $M_\sigma$, also denoted by $\Theta$.\\
\indent It remains only to show that $\Theta$ implements the conjugacy between $\rho|_{M_{\rho}}$ and $\sigma|_{M_\sigma}$.
Indeed,  for any positive integer $n$, and for any $x\in \rho^n(M)'\cap M$, we have\begin{eqnarray*}
\sigma\circ \Theta (x)&=&\sigma\circ\op{Ad}(u_n)(x)\\&=&\op{Ad}(u)\circ\rho\circ\op{Ad}(u_n)(x)\\&=&\op{Ad}(u_{n+1})\circ\rho(x)\\&=&\Theta\circ\rho(x),
\end{eqnarray*}where the last equality follows from (\ref{rel1}) and (\ref{rel2}). It then follows, by norm continuity, that $\sigma\circ\Theta(x)=\Theta\circ\rho(x)$ for all $x\in M_\rho$, so $\rho\mid_{{M}\sb{\rho}}$ and $\sigma\mid_{{M}\sb{\sigma}}$ are conjugate endomorphisms. \\
(II) For proving the general case, let $\sigma':=\theta\circ\sigma\circ\theta^{-1}\in\E$. Then by (I), the endomorphisms $\sigma'|_{M_{\sigma'}}$ and $\rho|_{M_{\rho}}$ are conjugate. But $\theta|_{M_{\sigma}}$  is also a *-isomorphism from $M_{\sigma}$ onto $M_{\sigma'}$ which implements the conjugacy  between $\sigma\mid_{{M}\sb{\sigma}}$ and $\sigma'\mid_{{M}\sb{\sigma'}}$. Thus $\rho\mid_{{M}\sb{\rho}}$and $\sigma\mid_{{M}\sb{\sigma}}$ are conjugate endomorphisms.\end{proof}
\section{Consequences}We infer immediately from Theorem \ref{cut1}, that in order to find classes of cocycle-conjugate endomorphisms of a von Neumann algebra $M$ that are also conjugate, we should look at those endomorphisms $\rho$ of $M$ that satisfy $$M_\rho=M.$$ An endomorphism $\rho\in\op{End}(M)$ that satisfies this relation is said to have the generating property (see also \cite{Ch94}). We note that the class of endomorphisms having the generating property is related to the class of Powers shift endomorphisms \cite{P88}, i.e., those endomorphisms $\rho$ of $M$ that satisfy the range property $\bigcap_n \rho^n(M)=\mathbb{C}1$:
\begin{lemma}\label{cut2}If $M$ is a factor, then any endomorphism $\rho$ of $M$ that has the generating property is a Powers shift endomorphism. If $M$ is a type $\rm{I}\sb{\infty}$ factor, then any Powers shift endomorphism has the generating property.
\end{lemma}
\begin{proof}If $\rho\in\op{End}(M)$, then it is easily seen that $$M_\rho\subset \left(\bigcap _{n\geq 0}\rho^n(M)\right)'\bigcap M.$$
Therefore, if $M$ is a factor and $M_\rho=M$, then $$
\bigcap _{n\geq 0}\rho^n(M)\subset M'\cap M=\mathbb{C}1,$$
so $\rho$ is a Powers shift.\\\indent If $M$ is a type $\rm{I}\sb{\infty}$ factor, and if $\rho$ is a Powers shift endomorphism of $M$, then, by using Corollary 2.2 of \cite{P88}, we have $$M=
\bigvee_{n\geq 0}\rho^n(\rho(M)'\cap M)),$$ the von Neumann algebra generated by $\{\rho^n(\rho(M)'\cap M)\,|\,n\geq 0\}.$ Moreover, since $$\rho^n(\rho(M)'\cap M)\subset \rho^{n+1}(M)'\cap M,\;\;\mbox{for all}\;n\geq 0,$$ we deduce that $\rho$ has the generating property. 
\end{proof}
\begin{remark}Examples of Powers shift endomorphisms that do not have the generating property were constructed explicitly in \cite{By88}, \cite{By90}.

\end{remark}
Let now $M$ be a type $\rm{I}\sb{\infty}$ factor. In \cite{A89}, W. Arveson has shown that any endomorphism $\rho$ of $M$ has the form $$
\rho(x)=\sum_{i=1}^ku_ixu_i^*, \;\;x\in M,$$ where $\{u_i\}_{i=\overline{1,k}}$ is a family of isometries in $M$ that satisfy the Cuntz relations \cite{C77}, and $k$ is the (Powers) index of $\rho$ \cite{P88}, i.e. that $k$ such that $\rho(M)'\cap M$ is a factor of type $\rm{I}\sb{k}.$\\\indent This index is a complete cocycle conjugacy invariant \cite[Theorem 2.4]{P88}, \cite[Proposition 2.3]{L93}. In fact, if $\rho(x)=\sum_{i=1}^ku_ixu_i^*$ and $\sigma(x)=\sum_{i=1}^kv_ixv_i^*$ are endomorphisms acting on $M$, then $$u=\sum_{i=1}^kv_iu^*_i$$ defines a unitary of $M$ such that $\sigma=\op{Ad}(u)\circ\rho$. \\
\indent The problem of finding conjugacy invariants is more difficult. For example, it was shown in \cite{B96} that it is impossible to find a smooth labeling of all the conjugacy classes of Powers shifts.
For this class of endomorphisms, the only known criterion for conjugacy is that found  by R. Powers \cite[Theorem 2.3]{P88}:  if $\rho$ and $\sigma$ are Power shifts with the same index, and if there exist a $\rho$-invariant pure normal state of $M$ and a $\sigma$-invariant pure normal state of $M$, then $\rho$ and $\sigma $ are conjugate endomorphisms. \\\indent Theorem \ref{cut1} and Lemma \ref{cut2} allow us to formulate a new criterion of conjugacy:  
\begin{corollary}Let $M$ be a type $\rm{I}\sb{\infty}$ factor, and let $\rho(x)=\sum_{i=1}^ku_ixu_i^*,\;\;\sigma(x)=\sum_{i=1}^kv_ixv_i^*,$ be Powers shift endomorphisms of index $k$ acting on $M$. Suppose that there exists a $\rho$-invariant faithful normal state $\vp$ of $M$ such that $$\sum_{i=1}^kv_iu^*_i\in M^\vp.$$Then $\rho$ and $\sigma$ are conjugate endomorphisms.

\end{corollary}
Let now $M$ be a von Neumann algebra. If the state $\vp$ in Theorem \ref{cut1} is tracial, then the condition $u\in M^\vp$ is superfluous. Moreover, if $M$ is a type $\rm{II}\sb{1}$ factor, and if the state $\vp$ is the unique faithful normal finite trace on $M$, then all conditions in Theorem \ref{cut1} that involve this state $\vp$ are automatically satisfied. Thus, we obtain the following consequence of Theorem \ref{cut1}, which generalizes Corollary 1.3 of \cite{By90}:
\begin{corollary}
Let $M$ be a von Neumann algebra, and let $\rho$, $\sigma$ be endomorphisms of $M$ having the generating property. Assume that there exists a $\rho$-invariant faithful normal tracial state of $M$. Then $\rho$ and $\sigma$ are conjugate endomorphisms if and only if they are cocycle-conjugate endomorphisms.
\end{corollary}

\bibliographystyle{amsalpha}

\end{document}